\documentclass[12pt]{article}
\usepackage{amsfonts}
\usepackage{epsfig}

\title{Rectangle Sweepouts and Coincidences}
\author{Richard Evan Schwartz \thanks{\hskip 5 pt Supported by 
N.S.F. Research Grant DMS-1204471}}

\newtheorem{theorem}{Theorem}[section]

\newtheorem{lemma}[theorem]{Lemma}

\newtheorem{corollary}[theorem]{Corollary}

\def\startproof{{\bf {\medskip}{\noindent}Proof: }}

\def\endproof{$\spadesuit$  \newline}

\def\R{\mbox{\boldmath{$R$}}}%
\def\Z{\mbox{\boldmath{$Z$}}}%

\begin{document}

\maketitle
\begin{abstract}
We prove an integral formula for continuous
paths of rectangles inscribed in a piecewise
smooth loop.  We then use this integral
formula to show that (with a very mild
genericity hypothesis) the number of
rectangle coincidences, informally described
as the number of inscribed
rectangles minus the number of isometry
classes of inscribed rectangles, grows
linearly with the number of positively oriented extremal
chords -- a.k.a. diameters -- in a polygon.
\end{abstract}

\section{Introduction}

A {\it Jordan loop\/} is the image of a circle under a
continuous injective map into the plane.
Toeplitz conjectured in 1911 that every Jordan loop contains $4$
points which are the vertices of a square.  This is sometimes called
the {\it Square Peg Problem\/}.  
For historical details and a long bibliography, we refer the
reader to the excellent survey article [{\bf M\/}] by B. Matschke,
written in 2014, and also Chapter 5 of I. Pak's online book
[{\bf P\/}]. 

Some interesting work on problems related to
the Square Peg Problem has been done very recently.
The paper of C. Hugelmeyer [{\bf H\/}] shows 
that a smooth Jordan loop
always has an inscribed rectangle of
aspect ratio $\sqrt 3$. The paper
[{\bf AA\/}] proves that any cyclic quadrilateral
can (up to similarity) be inscribed in any
convex smooth curve.  The
paper  [{\bf ACFSST\/}] proves, among other things,
that a dense set
of points on an embedded loop in space are
vertices of a (possibly degenerate) inscribed
parallelogram.

Say that a rectangle $R$ {\it graces\/} a
Jordan loop $\gamma$ if the vertices of
$R$ lie in $\gamma$ and if the 
cyclic ordering on the vertices induced by
$R$ coincides with the cyclic ordering
induced by $\gamma$.
Let $G(\gamma)$ denote the space of 
labeled gracing rectangles. 
In [{\bf S1\/}] we prove the following
result.

\begin{theorem} 
\label{threepoint}
Let $\gamma$ be a
Jordan loop.  Then $G(\gamma)$
contains a connected set $S$ such that
all but at most $4$ vertices of
$\gamma$ are vertices of members of $S$.
\end{theorem}

We have a more precise characterization of
the possibilities for $S$ in [{\bf S1\/}].
We proved Theorem
\ref{threepoint} by taking a limit of
a result for polygons.  We now
describe this result.

Given a polygon $P$, we say that a
chord $d$ of $P$ is a {\it diameter\/} 
if $d$ if the two perpendiculars to $d$
based at $\partial P$ do not
locally separate $\partial P$ into two arcs.
Each diameter can be positively oriented
or negatively oriented, but not both.
To explain the condition, we rotate
the picture so that $d$ is vertical.
The endpoints of $d$ divide $P$ into
two arcs $P_1$ and $P_2$.  Given the
non-separating condition associated
to a chord, we can say whether $P_1$
locally lies to the left or right of $P_2$
in a neighborhood of each endpoint of
$d$.  We call $d$ {\it positively oriented\/}
if the left/right answer is the same
at both endpoints. That is, either
$P_1$ locally lies to the left at
both endpoints or $P_1$ locally
lies to the right at both endpoints.
Figure 1 some examples of positive diameters.

\begin{center}
\resizebox{!}{2.5in}{\includegraphics{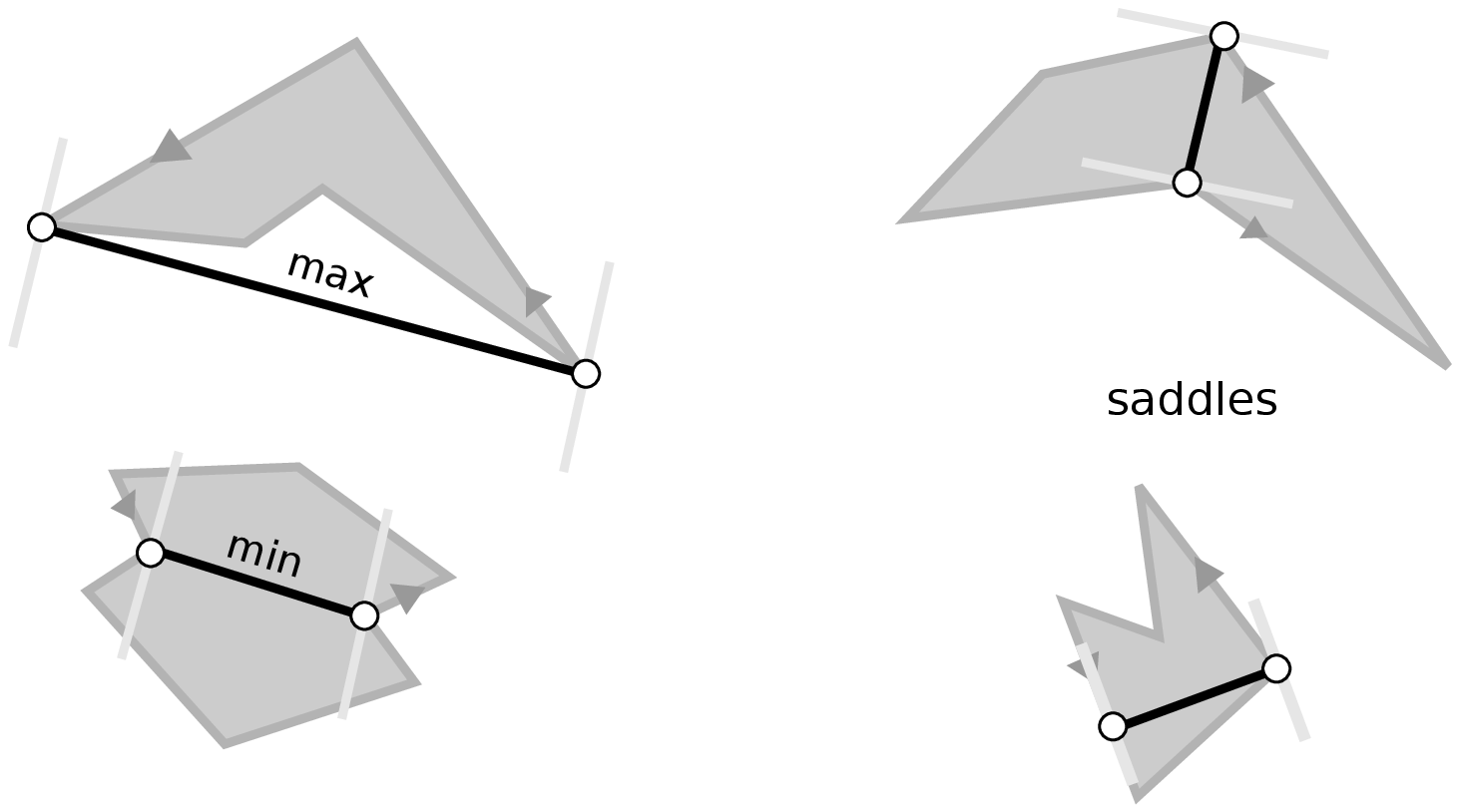}}
\newline
{\bf Figure 1:\/} Some positive diameters of polygons.
\end{center}

With respect to the distance function on
$P$, a diameter can be a minimum, a maximum,
or neither.  We call the third kind
{\it saddles\/}. 
Let $\Delta_+(P)$ denote the number of positively oriented
diameters of $P$.

Let $\Pi_N$ denote
the space of embedded $N$-gons.  The set $\Pi_N$
is naturally an open subset of
$(\R^2)^N$ and as such inherits the
structure of a smooth manifold.
We call a subset $\Pi_N^* \subset \Pi_N$ {\it fat\/}
if $\Pi_N-\Pi_N^*$ is a finite union of positive
codimension submanifolds. In particular,
a fat set is open and has full measure. 

\begin{theorem}
\label{polygon}
There exists a fat subset $\Pi_N^* \subset \Pi_N$
with the following property.  For every
$N$-gon $P \in \Pi_N^*$ the space
$\Gamma(P)$ is a piecewise-smooth $1$-manifold.
Each arc component of $\Gamma(P)$ connects
two positive diameters of $P$, and every positive diameter
arises as the end of $4$ arc components.
of $\Gamma(P)$.  In particular, there are
$2\Delta_+(P)$ arc components of $\Gamma(P)$.
\end{theorem}
The reason that there are $4$ arc components
connecting every pair of positive diameters that is
that we are considering cyclically labeled
rectangles.  Each of the $4$ components
is obtained from each other one by
cyclically relabeling.

Now we describe the results we prove in this paper.
Given a rectangle $R$, we let
$X(R)$ and $Y(R)$ respectively denote the
lengths of the first and second sides of $R$.
For any continuous path of rectangles
in $\Gamma(P)$ which is either a closed loop
or which connects two diameters of $P$, we define
the {\it shape curve\/}
$Z(\alpha)$.  This curve is given by
\begin{equation}
Z(\alpha,t)=(X(R_t),Y(R_t)).
\end{equation}
Here $t \to R_t$ is a parametrization of
$\alpha$.

When $\alpha$ is a closed loop,
$Z(\alpha)$ is a closed loop as well.
When $\alpha$ is an arc component,
$Z(\alpha)$ is an arc, not necessarily
embedded, that
starts and ends on the coordinate axes.
Figure 2 shows two of the possibilities.

\begin{center}
\resizebox{!}{1.5in}{\includegraphics{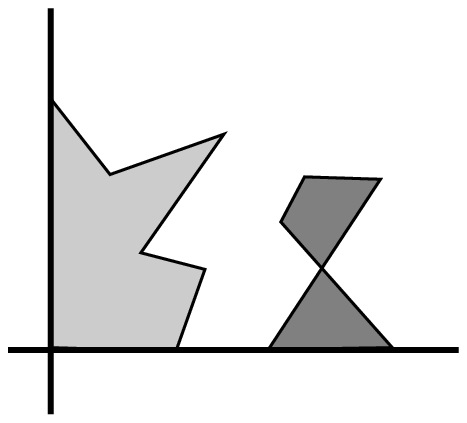}}
\newline
{\bf Figure 2:\/} Shape curves associated to 
hyperbolic and null arcs.
\end{center}

In the first case, one endpoint of
$\alpha$ lies on the $X$-axis and
the second endpoint lies on the $Y$-axis.
As in [{\bf S1\/}] we call such arcs
{\it hyperbolic arcs\/}.
In the other cases, both ends lie on the
same axis.  We call such components
{\it null arcs\/}.
In the arc cases, we
augment $Z(\alpha)$ by adjoining
the relevant parts of
the coordinate axes so as to create a
closed loop.  We have
shaded in the regions bounded by these
closed loops.  We call this augmented
loop the {\it shape loop\/} associated
to $\alpha$ and give it the same name.

In [{\bf S2\/}] we found a kind of integral
formula associated to the shape loop, though
we stated it in a different context.  
This invariant is quite similar to the
integral invariant in [{\bf Ta\/}], though
we use it in a different context.  (In 
\S \ref{squeeze} we give a sample result
from [{\bf S2\/}].)  Here we
adapt the invariant to the present situation
and prove the following theorem.

\begin{theorem}
\label{sweep}
Let $P$ be any piecewise smooth Jordan loop.
Let $\alpha$ be a piecewise smooth path in
$\Gamma(P)$.  If $\alpha$ is a hyperbolic
arc then the signed area of the region bounded by
$Z(\alpha)$ equals (up to sign) the area of
the region bounded by $P$.  If
$\alpha$ is either a null arc or a closed
loop, then the signed area of the region
bounded by $Z(\alpha)$ is $0$.
\end{theorem}

Theorem \ref{sweep} says something about the number
of coincidences that appear amongst the inscribed
rectangles.  We will give an example which
explains the connection. Since the shape loop associated to a
null component bounds a region of area $0$, the
shape curve must have a self-intersection. This
self-intersection corresponds to a pair of
isometric rectangles inscribed in the polygon.
Now we formulate a general result.
We call two labeled rectangles {\it really distinct\/}
if their unlabeled versions are also distinct.
Thus, two relabelings of the same rectangle are
not really distinct.

We define the multiplicity of the pair
$(X,Y)$ as follows.
\begin{itemize} 
\item $\mu(X,Y)=n-1$ if there are
$n>1$ really distinct
labeled rectangles $R_1,...,R_n$ inscribed in $P$ such
$X(R_j)=X$ and $Y(R_j)=Y$ for all $j=1,...,n$.
We also allow $n=\infty$,
\item $\mu(X,Y)=0$ if there are
$0$ or $1$ such rectangles.
\end{itemize}
We define
\begin{equation}
\label{coincidence}
M(P)=\sum \mu(X,Y),
\end{equation}
where the sum is taken over all pairs $(X,Y)$. Typically
this is a sum with finitely many finite nonzero terms.
There is a more natural (but somewhat informal)
way to think about $M(P)$.  Suppose that we
color all the points in $\Gamma(P)$ according to
the isometry class of rectangles they represent.
Then $M(P)$ is the number of points minus the
number of colors.

\begin{theorem}
\label{main}
For each $P \in \Pi_N^*$ we have
$M(P) \geq 2(\Delta_+(P)-2)$.
\end{theorem}
When $P$ is an obtuse triangle we have
$M(P)=0$ and $\Delta_+(P)=2$, so the result
is sharp in a trivial way.

Some version of Theorem \ref{main} is 
true for an arbitrary polygon, but here
we place a mild constraint so as to
make the proof easier.  Let $P$ be a polygon.
We call a diameter $S$ of $P$
{\it tricky\/} if the endpoints
of $S$ are vertices of $P$ and 
if at least one of the edges
of $P$ incident to $S$ is 
perpendicular to $S$.

\begin{theorem}
\label{main2}
If $P$ has no tricky diameters,
$M(P) \geq \frac{1}{16}(\Delta_+(P)-2)$.
\end{theorem}

The rest of the paper is devoted to proving the
results above.

\newpage

\section{The Integral Formula}

\subsection{The Differential Version}

Let $J$ be a piecewise smooth Jordan loop and
let $R$ be a labeled rectangle that graces $J$.
For each $j=1,2,3,4$ we let
$A_j$ denote the signed area of the region $R_j^*$ bounded
by the segment $\overline{R_{j}R_{j+1}}$ and
the arc of $J$ that connects $R_j$ to
$R_{j+1}$ and is between these two points in
the counterclockwise order.  Figure 3 shows a
simple example.  The signs are taken so that the signed areas
are positive in the convex case, and then in 
general we define the signs so that the signed
areas vary continuously.

\begin{center}
\resizebox{!}{1.9in}{\includegraphics{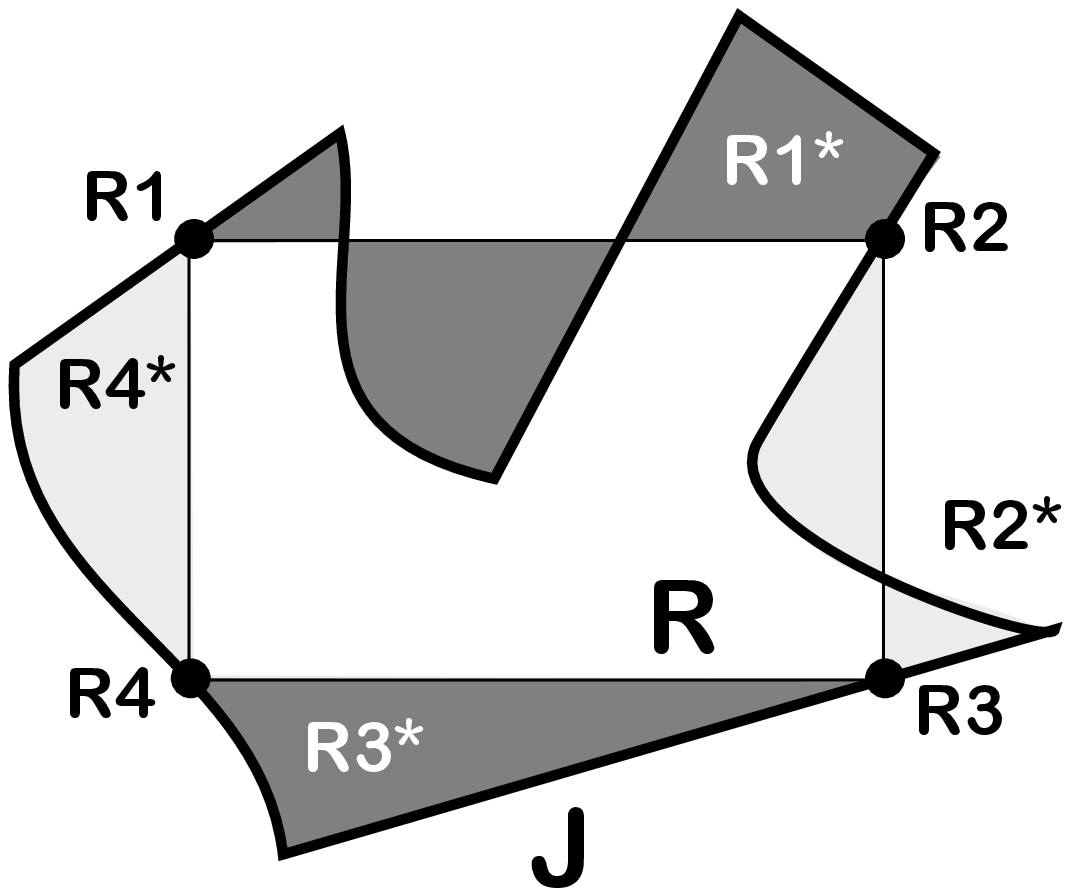}}
\newline
{\bf Figure 3:\/} The curve $J$, the rectangle $R$
and the regions $R_j^*$ for $j=1,2,3,4$.
\end{center}

Assuming that $J$ is fixed, we introduce the
quantity
\begin{equation}
A(R)=(A_1+A_3)-(A_2+A_4).
\end{equation}
We also have the point $(X,Y) \in \R^2$, where
\begin{equation}
X={\rm length\/}(\overline{R_1R_2}), \hskip 15 pt
Y={\rm length\/}(\overline{R_2R_3}), \hskip 15 pt
\end{equation}

Assuming that we have a piecewise smooth
path $t \to R_t$ of rectangles gracing $J$,
we have the two quantities
\begin{equation}
A_t=A(R_t), \hskip 30 pt
(X_t,Y_t)=(X(R_t),Y(R_t)).
\end{equation}

If $t$ is a point of differentiability, we may
take derivatives of all these quantities.
Here is the main formula.
\begin{equation}
\frac{dA}{dt}=Y \frac{dX}{dt} - X \frac{dY}{dt}.
\end{equation}

It suffices to prove this result for $t=0$.
This formula is rotation invariant, so for
the purposes of derivation, we rotate the
picture so that the first side of $R_0$ is
contained in a horizontal line, as shown
in Figures 3 and 4.  When we differentiate,
we evaluate all derivatives at $t=0$.

We write
\begin{equation}
\frac{dR_j}{dt}=(V_j,W_j).
\end{equation}

Up to second order, the region $R_1^*(t)$ 
is obtained by adding a small quadrilateral
with base $X_0$ and adjacent sides
parallel to $t(V_1,W_1)$ and $t(V_2,W_2)$.
Up to second order, the area of this
quadrilateral is 
$$\frac{X(W_1+W_2)}{2}.$$

\begin{center}
\resizebox{!}{2in}{\includegraphics{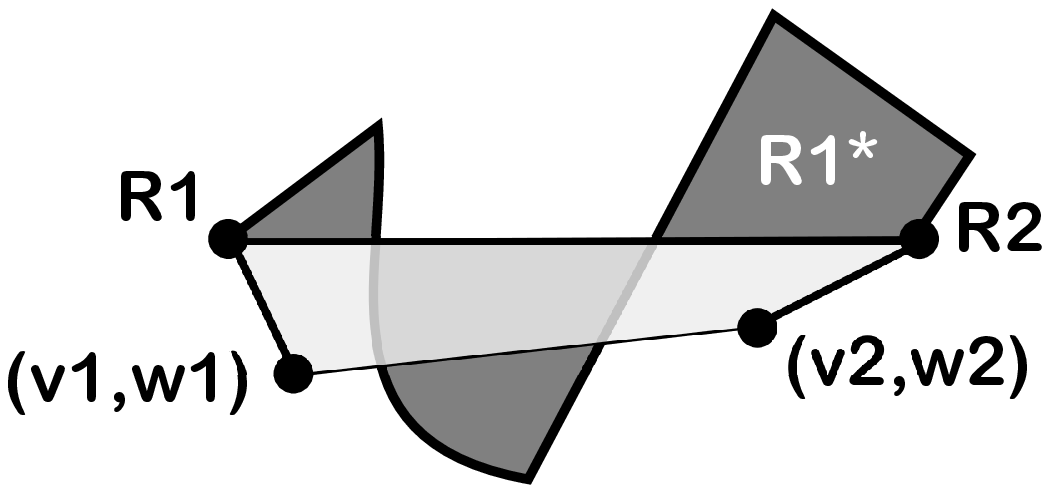}}
\newline
{\bf Figure 4:\/} The change in area.
\end{center}

From this equation, we conclude that
\begin{equation}
\frac{dA_1}{dt}=-\frac{X(W_1+W_2)}{2}.
\end{equation}
We get the negative sign because the area
of the region increases when $W_1$ and $W_2$
are negative.
A similar derivation gives
\begin{equation}
\frac{dA_3}{dt}=+\frac{X(W_3+W_4)}{2}.
\end{equation}
Adding these together gives
$$
\frac{dA_1}{dt}+\frac{dA_3}{dt}=
X \times \bigg[\frac{W_3-W_1}{2}\bigg] +
X \times \bigg[\frac{W_4-W_2}{2}\bigg]=$$
\begin{equation}
\label{term1}
-X \times \bigg[\frac{1}{2}\frac{dY}{dt}\bigg]+
-X \times \bigg[\frac{1}{2}\frac{dY}{dt}\bigg]=
-X \frac{dY}{dt}.
\end{equation}

A similar derivation gives
\begin{equation}
\frac{dA_2}{dt}=-\frac{X(V_2+V_3)}{2}, \hskip 30 pt
\frac{dA_4}{dt}=+\frac{X(V_4+V_1)}{2}.
\end{equation}
Adding these together gives
\begin{equation}
\label{term2}
\frac{dA_2}{dt}+\frac{dA_4}{dt}=
-Y \frac{dX}{dt}.
\end{equation}
Subtracting Equation \ref{term2} from
Equation \ref{term1} gives
\begin{equation}
\label{diff}
\frac{dA}{dt}=-X \frac{dY}{dt}+Y \frac{dX}{dt},
\end{equation}
as claimed.

\subsection{The Integral Version}

Let $\omega=-XdY+YdX$.  Here we think
of $\omega$ as a $1$-form.  Suppose that
we have parameterized our curve
of rectangles so that the parameter
$t$ runs from $0$ to $1$.  Integrating
Equation \ref{diff} over the piecewise
smooth path, we see that
\begin{equation}
\label{int}
A_1-A_0=\int_Z \omega.
\end{equation}
Here $Z$ is the shape curve associated to
the path of rectangles.  We can interpret
this integral geometrically.  Letting $O=(0,0)$,
consider the
closed loop
\begin{equation}
Z'=\overline{O, Z_0} \cup Z \cup \overline {Z_1,O}.
\end{equation}
Since $\omega$ vanishes on vectors of
the form $(h,h)$, we see that
\begin{equation}
A_1-A_0=\int_Z\omega=\int_{Z'} \omega= -\int \int_{\Omega} 2dxdy = -2\ {\rm area\/}(\Omega).
\end{equation}
Here $\Omega$ is the region bounded by $Z'$.
The last line of the equation refers to the signed area of $\Omega$.
\newline

\noindent
{\bf Proof of Theorem \ref{sweep}:\/}
Suppose first that $\alpha$ is a piecewise smooth loop rectangles
which grace the Jordan curve $J$.  Then the curve $Z$ is already
a closed loop, and the signed area of the region bounded by
$Z$ is the same as the signed area bounded by $Z'$.  Since
$A_1=A_0$ in this case, we see that $Z$ bounds a region of
signed area $0$.

If $\alpha$ is a null arc, then $R_0$ and $R_1$ both have
the same aspect ratio, either $0$ or $\infty$.  In either
case, we have $A_0=A_1$. The common value is, up to sign,
the area of the region bounded by $J$.  In this case,
$Z$ starts and stops on one of the coordinate axes, and
the region bounded by $Z$ has the same area as the 
shape loop, $Z \cup \overline{Z_0Z_1}$.  So, in this
case we also see that the shape loop bounds a region
of area $0$.

If $\alpha$ is a hyperbolic arc, then 
$A_0=-A_1$ and both quantities up to sign
equal the area of the region bounded by $J$.
At the same time $Z'$ is precisely the
shape loop in this case.   So, we see
that twice the area of the region bounded
by $J$ equals twice the area of the region
bounded by $Z$, up to sign.  Cancelling
the factor of $2$ gives the desired result.
\endproof

\subsection{Generic Coincidences}
\label{generic}

In this section we prove
Theorem \ref{main}.
Suppose that $P$ is an $N$-gon that satisfies
the conclusions of Theorem \ref{polygon}.
This happens if $P \in \Pi_N^*$, but it might
happen more generally.  In any case,
the space $\Gamma(P)$ of gracing
rectangles has $2\Delta(P)$ arc
components.  There is a $\Z/4$ action
on $\Gamma(P)$ and this action 
freely permutes the arc components
of $\Gamma(P)$.  

We let $\delta=\Delta/2$
and we let $\alpha_1,...,\alpha_{\delta}$
denote a complete set of representatives
of these arc components modulo the 
$\Z/4$ action.  It suffices to show that
the sum in Equation \ref{coincidence}
is at least $\delta-1$ when we
restrict our attention to the components
just listed.

Consider those arcs on our list which
are null arcs.  The shape
loops associated to each of these arcs
bound regions of area $0$ and hence
the corresponding loop has a double point.
Each double point corresponds to
a distinct pair that adds $1$ to the
total count for $M(J)$.  The remaining
rectangle coincidences involve rectangles
not associated to these arcs or to
their images under the $\Z/4$ action.

Now consider those arcs on our
list which are hyperbolic arcs
whose shape loops are not
embedded.  In exactly the same way
as above, each of these arcs
contributes $1$ to the count for
$M(J)$ and the rectangle pairs
involved are distinct from the
ones we have already considered.
Again, the remaining rectangle
coincidences involve rectangles
not associated to these arcs
or to their images under the $\Z/4$ action.
\newline
\newline
{\bf Remark:\/} 
Before we move on to the last
case, we mention that the
count above might be an under-approximation,
even in case there is just one double
point per shape loop considered.
Consider the simple situation
where there are just $2$ null
arcs.  It might happen that the
rectangle pairs corresponding
to these $2$ arcs are congruent
to each other.  This would give
us a $4$ congruent gracing
rectangles and would contribute
$3$ rather than $2$ to the
total count.
\newline

Finally, consider the $d$
hyperbolic arcs on our list which
have embedded shape loops.  If
$\alpha_1$ and $\alpha_2$ are
two such arcs, then
$Z(\alpha_1)$ and $Z(\alpha_2)$ are
two closed loops which bound the
same area. If these loops did not
intersect in the positive quadrant,
then either the region bounded
by $Z(\alpha_1)$ would strictly
contain the region bounded by
$Z(\alpha_2)$ or the reverse.
This contradicts the fact that
these two regions have the same
area.  Hence $Z(\alpha_1)$ and
$Z(\alpha_2)$ intersect in
the positive quadrant, and the
intersection point corresponds
to a coincidence involving a
rectangle associated to
$\alpha_1$ and a rectangle
associated to $\alpha_2$. 
Call this the {\it intersection property\/}.

We label so that
$\alpha_1,...,\alpha_d$ are the
hyperbolic arcs having embedded
shape loops.
We argue by induction that
these $d$ arcs contribute at
least $d-1$ to the count for
$M(J)$.  If $d=1$ then there
is nothing to prove.
By induction, rectangle
coincidences associated to the arcs
$\alpha_1,...,\alpha_{d-1}$
contribute $d-2$ to the count
for $M(J)$.

By the intersection property,
$\alpha_d$ intersects each
of the other arcs, and
$\Gamma(J)$ is a manifold, there
is at least one new rectangle
involved in our count, namely
one that corresponds to a point
on $Z(\alpha_d)$ that is also
on some of the shape loop.
The corresponding rectangle
adds $1$ to the count in
Equation \ref{coincidence},
one way or another.  So,
all in all, we add
$d-1$ to the count for
$M(J)$ by considering
the rectangle coincidences
associated to $\alpha_1,...,\alpha_d$.
This proves what we want.

\subsection{A Non-Squeezing Result}
\label{squeeze}

Here we explain how the invariant above
implies one of our main results in
[{\bf S2\/}].  Really, it is the
same proof.  The material
in this section plays no role in the
rest of the paper.

Suppose that $\gamma_1$ and $\gamma_2$ are
$2$ piecewise smooth curves which are
disjoint.  Suppose also that at each
end, $\gamma_j$ coincides with a
line segment. Finally suppose that
these line segments are parallel
at each end, so to speak. Figure 5 shows
what we mean.

\begin{center}
\resizebox{!}{1.5in}{\includegraphics{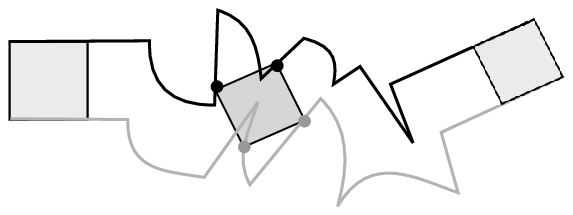}}
\newline
{\bf Figure 5:\/} Sliding a square along a track.
\end{center}

Suppose that we have a piecewise smooth family
of rectangles, all having the same
aspect ratio, that starts at one end,
finishes at the other, and remains inscribed
in $\gamma_1 \cup \gamma_2$ the whole time.
We imagine $\gamma_1 \cup \gamma_2$ as being
a kind of track that the rectangle slides
along (changing its size and orientation
along the way).
Figure $5$ shows an example in which case the
rectangle is a square.  In Figure 5
we show the starting rectangle $R_0$, the
ending rectangle $R_1$, and some $R_t$ for
$t \in (0,1)$.  This is just a hypothetical
example.

We can complete the union $\gamma_1 \cup \gamma_2$ to
a piecewise smooth Jordan loop by extending
the ends of one or both of these curves,
if necessary, and then
dropping perpendiculars. Let $\Omega$ be
the region bounded by this loop.  The shape curve
associated to our path lies on a line through
the origin, and our $1$-form $\omega$ vanishes
on such lines.  Referring to the invariant
above, we therefore have  $A(R_0)=A(R_1)$.  But, after
suitably labeling the rectangles in our family, we have
$$A(R_j)={\rm area\/}(\Omega)-{\rm area\/}(R_j).$$
Hence $R_0$ and $R_1$ have the same area.
Since they also have the same aspect ratio,
they have the same side-lengths.
This is to say that the perpendicular
distance between the end of $\gamma_1$
and the end of $\gamma_2$ is the same
at either end.  This is a kind of
non-squeezing result.

In particular, our result shows that
Figure 5 depicts an impossible situation.
There is no way to slide a square
continuously through the shown ``track''
because the widths are different at the
$2$ ends.

\newpage

\section{The General Case}

\subsection{Rectangles Inscribed in Lines}
\label{conn}

The goal of this chapter is to prove
Theorem \ref{main2}.  We plan to
take a limit of the result in 
Theorem \ref{main}.

Let $E=(E_1,E_2,E_3,E_4)$ be  a
collection of $4$ line segments,
not necessarily distinct.
We say that a rectangle $R$ 
{\it graces\/} $E$ if the vertices
$R_1,R_2,R_3,R_4$
of $R$ go in cyclic order,
either clockwise or counterclockwise,
and $R_i \in E_i$ for all $i=1,2,3,4$.
We allow $R$ to be degenerate. 
Let $\Gamma(E) \subset (\R^2)^4$ denote the set of
rectangles gracing $E$. Note

We call a point $p \in \Gamma(E)$ 
{\it degenerate\/} if every neighborhood
of $p$ in $\Gamma_E$ contains points
corresponding to infinitely many
distinct but isometric rectangles.
We call $E$ {\it degenerate\/} if
there is some $p \in \Gamma(E)$ which
is degenerate.  

\begin{lemma}
Suppose that $E$ is nondegenerate. 
$\Gamma(E)$ is the intersection of
a conic section with a rectangular solid.
\end{lemma}

\startproof
Let $E=(E_1,E_2,E_3,E_4)$ be a
$4$-tuple of lines.  We rotate so
that none of the segments is vertical,
so that we may parameterize the
lines containing our segments
by their first coordinates.
Let $L_j$ be the line extending $E_j$.
We identify $\R^3$ with triples
$(x_1,x_2,x_3)$ where $p_j=(x_j,y_j) \in L_j$.
We let $p_4$ be such that
$p_1+p_3=p_2+p_4$.  In other words,
we choose $p_4$ to that
$(p_1,p_2,p_3,p_4)$ is a
parallelogram.

Let $\Gamma(L)$ denote the set of
rectangles gracing $L$.  We describe
the subset $\Gamma'(L) \subset \R^3$ corresponding
to $\Gamma(L)$.  The actual
set $\Gamma(L)$ is the image of
$\Gamma'(L)$ under a linear map from
$\R^3$ into $(\R^2)^4$.

The condition that $p_4 \in L_4$ is a linear
condition.  Therefore, the set
$(x_1,x_2,x_3) \in \R^3$ corresponding to
parallelograms inscribed in $L$ is a
hyperplane $\Pi$.   The condition that
our parallelogram is a rectangle is 
$(p_3-p_2) \cdot (p_1-p_2)=0.$
This condition defines a quadric
hypersurface $H$ in $\R^3$. 
The intersection $\Gamma'(L)=\Pi \cap H$
corresponds to the
inscribed rectangles. 

$\Pi \cap H$ is either a plane
or a conic section.  In the
former case, $E$ is degenerate.
In the latter case, every point
$\Pi \cap H$ is either an
analytic curve or two crossing
lines.  Since $\Gamma(L)$ is
the image of $\Gamma'(L)$ under
a linear map, the set
$\Gamma(L)$ is also a conic section.

Let
$[E]=E_1 \times E_2 \times E_3 \times E_4.$
The $[E]$ is a rectangular solid. We have
$\Gamma(E)=\Gamma(L) \cap [E]$.
\endproof

\begin{lemma}
When $E$ is non-degenerate,
$\Gamma(E)$ has at most
$64=2^8$ connected components.
\end{lemma}

\startproof  
We use the notation from the previous
lemma. Note $[E]$ is
bounded by $8$ hyperplanes and
a conic section either lies in a
hyperplane or intersects it at
most twice.  So, each boundary
component of $[E]$ cuts
$\Gamma(L)$ into at most $2$ components.
\endproof

We call a polygon $P$
{\it degenerate\/} if some
$4$-tuple of edges associated to
$P$ is degenerate.
Otherwise we call $P$ {\it non-degenerate\/}.

\begin{lemma}
Let $P$ be a non-degenerate polygon.
The space $\Gamma(E)$ is a graph
having analytic edges and degree
at most $32$.
\end{lemma}

\startproof
Each rectangle $R$ can grace
at most $16$ different
$4$-tuples of edges of $P$, because each
vertex can lie in at most $2$ segments.
Hence, each $p \in \Gamma(P)$ lies
in the intersection of at most
$16$ distinct $\Gamma(E)$.
Since $\Gamma(E)$ is the intersection
of a conic section with a rectangular
solid, $\Gamma(E)$ is a graph with
analytic edges and maximum degree $4$.
From what we have said above,
$\Gamma(P)$ is a graph with
analytic edges and maximum degree $64=16 \times 4$.

We can cut down by a factor of
$2$ as follows.  The only time a
point of $\Gamma(P)$ lies in more
than $8$ spaces $\Gamma(E)$ is when
$p$ corresponds to a gracing
rectangle whose every vertex is a
vertex of $P$.  In this case,
$p$ is a vertex of $[E]$ for each
$4$-tuple $E$ that the rectangle
graces. But then
$p$ has degree at most $2$ in each
$\Gamma(E)$.  So, this exceptional
case produces vertices of degree
at most $32$.
\endproof

\subsection{The Inscribing Sequence}

A generic polygon $P$ satisfies the conclusions
of Theorem \ref{main}.  For such polygons,
any $4$-tuple which supports a
gracing rectangle is nice.

We label the sides of
$P$ by $\{1,...,N\}$.
Let $\Omega$ denote the set of ordered
$4$-element subsets of $\{1,...,N\}$,
not necessarily distinct.
Consider some embedded arc $\alpha \subset
\Gamma(P)$ of inscribed rectangles.
$\alpha$ defines a finite sequence
$\Sigma$ of elements of $\Omega$.
We simply note which
edges of $P_n$ contain any given
rectangle and then we order the
elements of $\Omega$ we get.
We call $\Sigma$ the
{\it inscribing sequence\/}
for $\alpha$.

\begin{lemma}
\label{inscribing}
$\Sigma$ has length at most $\kappa N^4$.
\end{lemma}

\startproof
If $\Sigma$ had length longer
than this, then we could find
a single $4$-tuple $E$ of edges
such that the subset of
$\alpha$ supported by $E$
has at least $82$ components.
In other words the sequence
would have to return to the
$4$-tuple describing $E$
at least $82$ times. The
arcs of
$\Gamma(E)$ corresponding to
these returns are disconnected
from each other, because otherwise
$\alpha$ would be a loop rather
than an arc.  This contradiction
proves our claim.
\endproof

\subsection{Stable Diameters}

For the rest of the chapter, we
use the word {\it diameter\/} to
mean a positively oriented diameter,
in the sense discussed in the introduction.

Let $P$ be a polygon and let
$S$ be a diameter of $P$.  We
call $S$ {\it stable\/} if
\begin{itemize}
\item At least one endpoint of $S$
is a vertex of $P$.
\item If $v$ is an endpoint of $S$
and $e$ is an edge of $P$ incident
to $P$ at $v$, then $S$ and $e$ are
not perpendicular.
\end{itemize}

\begin{lemma}
Suppose that $P$ has no tricky diameters.
If $P$ has an unstable diameter, then
$P$ is non-degenerate.
\end{lemma}

\startproof
This is a case-by-case analysis.
Suppose first that $P$ has a diameter $S$ whose
endpoints are not vertices of $P$.
Then the endpoints of $S$ lie in
the interior of a pair of parallel
edges of $P$.  But then $P$ is degenerate.
Suppose that $P$ has a diameter $S$ having
one endpoint which is a vertex $v$ of $P$.
The other endpoint of $S$ lies in the 
interior of an edge $e'$ of $P$. By
definition $e'$ and $S$ are perpendicular. If
$S$ is not stable, then one of the edges
$e$ of $P$ is perpendicular to $S$ and
hence parallel to $e'$. But then we can
shift $S$ over a bit and produce a diameter
of $P$ whose endpoints lie in the interior
of $e$ and $e'$.   Again, $P$ is
degenerate.  The remaining unstable
diameters are (in the technical sense) tricky.
\endproof

In view of the preceding result, it suffices to
prove Theorem \ref{main2} under the assumption
that $P$ is non-degenerate and has all stable
diameters.

\subsection{Limits of Diameters}

Let $P$ be an $N$-gon with
stable diameters.
We can find a sequence
$\{P_n\}$ of generic
$N$-gons converging to $P$.  Each
$P_n$ satisfies the conclusions
of Theorem \ref{main}.  

\begin{lemma}
Let $D$ be a diameter of $P$.
The polygon $P_n$ has a diameter $D_n$
such that $\{D_n\}$ converges to $D$.
\end{lemma}

\startproof
Since $P$ only has stable
diameters, there are
just $2$ cases to consider.
Suppose first that $D$ connects 
two vertices $v$ and $w$ of $P$.
The polygon
$P_n$ has vertices $v_n$ and $w_n$ which
converge respectively to $v$ and $w$
as $n \to \infty$.
Let $D_n$ be the chord whose endpoints
are $v_n$ and $w_n$.  By construction,
$D_n$ converges to $D$ and
for large $n$ this chord is a diameter.

Suppose now that $D$ connects a
vertex $v$ to a point in the interior
of an edge $e$.
Let $v_n$ and $e_n$
be the corresponding vertex and
edge of $P_n$. Since $v_n \to v$
and since $e_n \to e$ we see that
eventually there is a chord
$D_n$ that has $v_n$ as one endpoint
and has the other endpoint perpendicular
to $e_n$.  By construction $D_n \to D$
and eventually $D_n$ is a diameter
of $P_n$.
\endproof

\begin{lemma}
If $\{D_n\}$ is a sequence
of diameters of $P_n$, then 
$\{D_n\}$ converges on a subsequence to a diameter of $P$.
\end{lemma}

\startproof
Given the
sequence $\{D_n\}$ we can pass to a
subsequence so that the endpoints of
these diameters converge.  The limiting
segment $D$, provided that it has nonzero
length, must be a diameter of $P$ because
the required condition is a closed condition.
We just have to see that the length
of $\{D_n\}$ does not shrink to $0$.
Note that $D_n$ is at least as long
as the shortest diameter of $P_n$.
Furthermore, there is a positive
lower bound to the length of any
edge of $P_n$, independent of $n$.
So, if the length of $D_n$ converges
to $0$, there are two non-adjacent
vertices of $D_n$ whose distance
converges to $0$.  This contradicts
the fact that $\{P_n\}$ converges
to the embedded polygon $P$.
\endproof

We think of a diameter as a subset of
$(\R^2)^2$, and in this way we can
talk about the distance between
two diameters of $P_n$.

\begin{lemma}
Suppose that $\{D_n\}$ and
$\{D_n'\}$ are two sequences
of diameters such that the
distance from $D_n$ to $D_n'$
converges to $0$ as $n \to \infty$.
Then $D_n=D_n'$ for $n$ sufficiently
large.
\end{lemma}

\startproof
Let $v_n$ and $w_n$ be the
endpoints of $D_n$ and
let $v_n'$ and $w_n'$ be
the endpoints of $D_n'$.
We label so that
$\|v_n-v_n'\|$ and
$\|w_n-w_n'\|$ both tend to $0$.
In all cases, we can re-order so
that $v_n$ is a vertex of $P_n$
and $v_n'$ is not.  In other
words, $v_n'$ lies in the interior
of an edge $e_n'$ of $P_n$.
Since $v_n'$ converges to $v_n$,
a vertex of $P_n$, the segment
$e_n'$ becomes perpendicular to
$D_n'$ in the limit. This contradicts
the fact that $P$ has only stable
diameters.
\endproof

\begin{corollary}
\label{stable}
For $n$ sufficiently large, there
is a bijection between the diameters
of $P_n$ and the diameters of $P$
such that each diameter of $P$ 
is match with a sequence of diameters
of $P_n$ which converges to $P$.
\end{corollary}

\startproof
This is an immediate consequence of the
preceding $3$ lemmas.
\endproof

We truncate our sequence of polygons so that
the last corollary holds for all $n$.
For each $n$, these diameters 
are paired together by the arc components
of the manifold $\Gamma(P_n)$.  We pass to
a further subsequence so that the
same pairs arise for each $n$.  This gives
us a well defined way to pair the
diameters of $D$.  We say that two
diameters of $D$ are
{\it partners\/} if and only if the
corresponding diameters of $D_n$ are
paired together.

\begin{lemma}
\label{connect}
Each pair of partner diameters in
$P$ is connected by a piecewise smooth
path in $J(P)$.
\end{lemma}

\startproof
Let $A$ and $B$ be two partner diameters
of $P$.  Let $A_n$ and $B_n$ be the
corresponding diameters of $P_n$.
Let $\alpha_n$ be the arc in
$\Gamma(P_n)$ which connects
$A_n$ and $B_n$.
To understand the convergence of
$\{\alpha_n\}$ we
work in the Hausdorff topology
on the set of compact subsets of
$(\R^2)^4$.  This ambient
space contains $\Gamma(J)$ for
any Jordan loop.

We consistently label the
sides of $P_n$ and $P$.
Let $\Sigma_n$ be the
inscribing sequence of
$\alpha_n$.  By Lemma
\ref{inscribing} there is
a uniform upper bound of
$\kappa N^4$ on the length
of $\Sigma_n$. Therefore, we
may pass to a subsequence
so that the inscribing
sequence associated to 
$\alpha_n$ is independent of $n$.
We write
$$\alpha_n=\alpha_{n1},...,\alpha_{nk},$$
where $\alpha_{nj}$ is the arc of
rectangles corresponding to the
$j$th element of the sequence in
$\Omega$.  Here $k$ is the length
of the inscribing sequence.

We pass to a subsequence so that
$\{\alpha_{nj}\}$ converges in
the Hausdorff topology to a
subset $\alpha_j \subset \alpha$.
The set $\alpha_j$ is connected
and contained in a subset of
$\Gamma(E)$, where $E$ is the
$4$-tuple of edges corresponding
to the $j$th element of $\Omega$.
From the discussion in \S \ref{conn}, we see that
$\alpha_j$ is a compact, connected
algebraic arc.  By construction
$\alpha_j$ and $\alpha_{j+1}$ share
one point common for all $j$.
This vertex is the limit of
the sequence $\{\alpha_{nj} \cap \alpha_{n,j+1}\}$.

The description above reveals $\alpha$
to be a piecewise smooth arc
connecting the two diameters $A$ and $B$.
\endproof

\subsection{The End of the Proof}

Let $P$ be a polygon.
We still assume that $P$ has
stable diameters, so that the
results from the previous section apply.
We know from Lemma \ref{connect}
that the diameters of $P$
are paired in some way, and
each pair is connected by some
piecewise smooth path of
gracing rectangles.  We can
erase any loops that these
paths have and thereby
assume that all these paths
are embedded.  Next, we can assume that
every $2$ arcs in the collection
intersect each other in at most
one point.  Otherwise, we can
do a splicing operation to decrease
the number of intersection points.
(See Figure 6 below.)
The splicing operation may change
the way that the diameters are
paired up, but this doesn't bother us.
Finally, we
can make our choice of connectors
invariant under the
$\Z/4$ re-labelling action.

As in the proof of
Theorem \ref{main} we let
$\delta=\Delta_+(P)/2$ and we
chose a collection
$\alpha_1,...,\alpha_{\delta}$ of
connecting arcs which has
one representative in each
orbit of the $\Z/4$ action.

Suppose that our collection of
paths contains two hyperbolic
arcs $\alpha_1$ and $\alpha_2$
that intersect.  Each path
connects a (degenerate)
rectangle of aspect ratio $0$
to a (degenerate) rectangle
of aspect ratio $\infty$.
By splicing the paths together
and then re-dividing them,
we produce $2$ new paths
$\beta_1$ and $\beta_2$ such
that each $\beta_j$ connects
two degenerate rectangles of
the same aspect ratio.
In other words, we can do
a cut-and-paste operation
at an intersection point
to replace the two hyperbolic
arcs by null arcs.  If
necessary, we can erase any
loops created in this process.
Figure 6 shows this operation.

\begin{center}
\resizebox{!}{1.2in}{\includegraphics{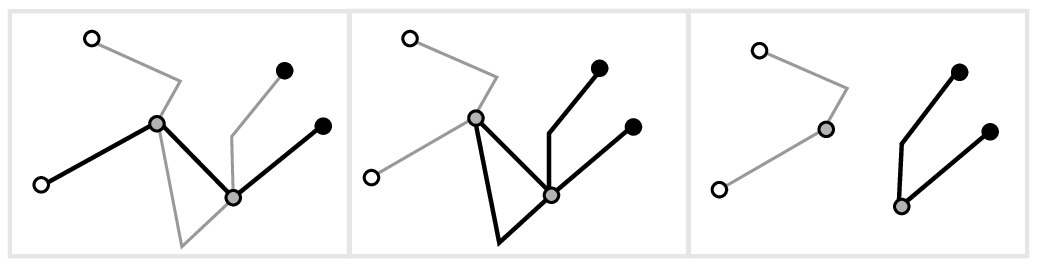}}
\newline
{\bf Figure 6:\/} The splicing operation.
\end{center}

Suppose first that there are $\delta/2$ arcs
in our collection that are hyperbolic
arcs.  Then this collection is
an embedded $1$-manifold contained in
$\Gamma(P)$.  Just using these arcs,
the same argument as in the
proof of Theorem \ref{main}
shows that
$$M(P) \geq \Delta_+(P)-2.$$
That is, we get the same answer
as in Theorem \ref{main} except
for the factor of $1/2$.

Now suppose that there are at least
$\delta/2$ null arcs.  For the rest
of the proof we just deal with these
null arcs.
Let $\Gamma_1(P)$ denote the
union of these null arcs.
We know that $\Gamma_1(P)$ is a
subset of $\Gamma(P)$ and
also a graph with algebraic
edges and maximim valence
at most $32$.
Let $\widehat \Gamma_1$ denote the
formal disjoint union of these
embedded null arcs.  The space
$\widehat \Gamma_1$ is a $1$-manifold,
just a union of arcs, and the 
``forgetful map'' $\phi: \widehat \Gamma_1 \to \Gamma_1$
is at most $16$ to $1$.

The same argument as in the proof
of Theorem \ref{generic} says that
there are $\delta$ distinct points
in $\widehat \Gamma_1$, two per
arc, corresponding to rectangle
coincidences.  Let $S$ be the
set of these points.  The image
$\phi(S)$ contains at least $\delta/16$
points.   For each of these points,
there is a second point corresponding
to an isometric rectangle. We know
this because the map
$\phi$ is injective on each null arc,
and each null arc contains $2$ 
points of $S$.  So, we can match
our $\delta/16$ points into
$\delta/32$ distinct pairs of
points, corresponding to pairs
of isometric but distinct
rectangles in $\Gamma(P)$.
This adds a count of $\delta/32$ to
$M(P)$. To make the comparison with
Theorem \ref{main} cleaner, we work
with $(\delta-1)/32$ instead.

In the case at hand, we get the same bound
as in Theorem \ref{main} except 
for the factor of $1/32$.
Going back to the count of labeled
rectangles, we have
$$M(P) \geq \frac{1}{16}(\Delta_+(P)-2).$$
This completes the proof of
Theorem \ref{main2}.

\newpage

\section{References}

[{\bf AA\/}] A. Akopyan and S Avvakumov, {\it Any cyclic quadrilateral can be
inscribed in any closed convex smooth curve.\/} 
arXiv: 1712.10205v1 (2017)
\newline
\newline
[{\bf ACFSST\/}] J. Aslam, S. Chen, F. Frick, S. Saloff-Coste,
L. Setiabrate, H. Thomas, {\it Splitting Loops and necklaces:
Variants of the Square Peg Problem\/},  arXiv 1806.02484 (2018)
\newline
\newline
[{\bf CH\/}] D. Hilbert and S. Cohn-Vossen, {\it Geometry and The Imagination\/},
\newline
Chelsea Publishing Company (American Math Society), 1990
\newline
\newline
[{\bf H\/}] C. Hugelmeyer, 
{\it Every Smooth Jordan Curve has an inscribed
rectangle with aspect ratio equal to $\sqrt 3$.\/}
arXiv 1803:07417 (2018)
\newline
\newline
[{\bf M\/}] B. Matschke, {\it A survey on the
Square Peg Problem\/},  Notices of the A.M.S.
{\bf Vol 61.4\/}, April 2014, pp 346-351.
\newline
\newline
[{\bf S1\/}] R. E. Schwartz, {\it A Trichotomy for Rectangles Inscribed in
Jordan Loops\/},  preprint, 2018
\newline
\newline
[{\bf S2\/}] R. E. Schwartz, {\it Four lines and a Rectangle\/},
preprint, 2018
\newline
\newline
[{\bf Ta\/}], T. Tao, {\it An integration approach
to the Toeplitz square peg conjecture\/}
\newline
Forum of Mathematics, Sigma, 5 (2017)
\newline
\newline
[{\bf W\/}] S. Wolfram, {\it The Mathematica Book\/}, 4th ed. Wolfram Media/Cambridge
University Press, Champaign/Cambridge (1999)

\end{document}